\input amstex
\documentstyle{amsppt}
\document
\magnification=1200
\NoBlackBoxes
\nologo
\pageheight{18cm}


\bigskip

\centerline{\bf $F$--MANIFOLDS WITH FLAT STRUCTURE}

\medskip

\centerline{\bf AND DUBROVIN'S DUALITY}

\bigskip 

\centerline{\bf Yuri I.~Manin}

\medskip

\centerline{\it MPIM, Bonn, Germany, and Northwestern University, Evanston, USA}

\bigskip

\medskip

{\bf Abstract.} This work continues the study of $F$--manifolds $(M,\circ )$,
first defined in [HeMa] and investigated in [He]. 
The notion of a compatible flat structure $\nabla$ is introduced,
and it is shown that many constructions known for Frobenius manifolds
do not in fact require invariant metrics and can be developed
for all such triples $(M,\circ ,\nabla ).$ In particular,
we extend and generalize recent Dubrovin's duality [Du2].

\bigskip

\centerline{\bf \S 1. Introduction.} 

\medskip

The notion of {\it Frobenius (super)manifold} $M$
axiomatized and thoroughly studied by B.~Dubrovin in the early nineties,
plays a central role in mirror symmetry, theory of unfolding
spaces of singularities, and quantum cohomology. 
A full Frobenius structure on $M$ consists of
the data $(\circ, e, g, E)$. Here $\circ :\Cal{T}_M\otimes_{\Cal{O}_M} \Cal{T}_M \to \Cal{T}_M$ is an associative and (super)commutative 
multiplication on
the tangent sheaf, so that $\Cal{T}_M$ becomes a sheaf of
(super)commutative $\Cal{O}_M$--algebras with identity $e\in \Cal{T}_M(M)$;
$g$ is a metric on $M$ (nondegenerate quadratic form $S^2(\Cal{T}_M)\to \Cal{O}_M$), and $E$ is an Euler vector field. These structures are
connected by various constraints and compatibility conditions,
spelled out in [Du1] and [Ma],
for example, $g$ must be flat and $\circ$--invariant.

\smallskip

Pretty soon it became clear  that various weaker versions of 
the Frobenius structure are interesting in themselves and also appear
naturally in different contexts. Here I focus on
the core notion of {\it $F$--manifold} introduced in [HeMa] and
further studied in [He]. This structure consists of
an associative and (super)commutative 
multiplication $\circ$  on
the tangent sheaf as above, constrained by the following
identity (1.2).  

\smallskip

Start with an expression
measuring the deviation
of the structure $(\Cal{T}_M, \circ , [,\,])$ from that
of a (sheaf of) Poisson algebra(s) on  $(\Cal{T}_M, \circ )$:
$$
P_X(Z,W):= [X,Z\circ W]-[X,Z]\circ W-(-1)^{XZ}Z\circ [X,W].
\eqno(1.1)
$$
Here $X,Y,Z,W$ are arbitrary local vector fields,
and a notation like $(-1)^{XZ}$ is a shorthand
for $(-1)^{\widetilde{X}\widetilde{Z}}$
where  $\widetilde{X}$ is the parity of $X$.

\smallskip 

Then we must have
$$
P_{X\circ Y}(Z,W)=X\circ P_Y(Z,W)+(-1)^{XY}Y\circ P_X(Z,W).
\eqno(1.2)
$$
If $\circ$ has an identity $e\in \Cal{T}_M$, we will
call $(M,\circ ,e)$ an $F$--manifold with identity.
\smallskip

$F$--manifolds keep reappearing in recent research,
although they are not always recognized as such.
Any Frobenius manifold stripped of $g$,$E$ and $e$,
becomes an $F$--manifold.
Solutions of the oriented
associativity equations with flat identity defined and studied
in [LoMa2], 5.3.1--5.3.2, are exactly $F$--manifolds with
compatible flat structure (see Definition 2.2 below).
Quantum $K$--theory produces $F$--manifolds 
with a flat invariant metric which are not quite
Frobenius because $e$ is not flat: see [Lee]. Dubrovin's 
almost Frobenius manifolds ([Du2], sec. 3, Definition 9)
are $F$--manifolds with nonflat identity as well.

\smallskip

In a very general context of $dg$--extended deformation
theory, S.~Merkulov in [Me1], [Me2], found out that extended
moduli spaces (e.~g. deformations of complex or symplectic structure)
often carry a natural $F$--structure, and produced a strong homotopy
version of the equation (1.2).
\smallskip

In this paper I introduce and study
$F$--manifolds {\it with a compatible flat structure.}
This notion turns out to share some of the deeper
properties of Frobenius manifolds, in particular, Dubrovin's
deformed connections, Dubrovin's duality, and several versions
of an operadic description. 

\smallskip

The paper is structured as follows. 
\smallskip

The main result
of \S 2 shows that there is an essential equivalence between
the notions of an $F$--manifold with a compatible flat
structure and that of a pencil of flat torsionless
connections. This fact then allows us to borrow
various techniques from [LoMa1] and
[LoMa2]. In \S 3 I generalize to this setup the
formalism of extended structure connections and show,
in particular, how Euler fields emerge as classifiers
of certain extended connections. Dubrovin's duality
is treated in \S 4: the approach via external vector bundles
with a pencil of flat connections makes Dubrovin's
construction more general and more transparent.
Finally, in  \S 5 I reproduce without proofs a
representation theoretic description of formal flat
manifolds from [LoMa1].

\bigskip

\centerline{\bf \S 2. $F$--manifolds with compatible flat structures } 

\medskip

 Below 
I broadly follow the
conventions of [Ma], Chapter I. 
``A manifold'' may mean an analytic supermanifold
over $\bold{C}$ (eventually over an exterior $\bold{C}$--algebra
carrying odd constants), or a germ of it, or a formal
completion of it along an embedded closed submanifold.
The sign rules are determined
by the convention that De Rham differentials and connections are odd.

\medskip

{\bf 2.1. Compatible flat structures.}  An (affine) flat structure on a manifold $M$ can be described
in three equivalent ways:

\smallskip

(i) A torsionless flat connection $\nabla_0:\,\Cal{T}_M\to
\Omega^1_M\otimes_{\Cal{O}_M}\Cal{T}_M $.

\smallskip

(ii) A local system $\Cal{T}_M^f\subset \Cal{T}_M$ of flat vector fields,
which forms a sheaf of supercommutative Lie algebras of rank $\roman{dim}\,M$
such that $\Cal{T}_M=\Cal{O}_M\otimes \Cal{T}_M^f$.

\smallskip

(iii) An atlas whose transition functions are affine linear.

\smallskip

The equivalence is established as follows: given $\nabla_0$, we put $\Cal{T}_M^f:= \roman{Ker}\,\nabla_0$; given a flat local map $(x^a)$,
we define local sections of $\Cal{T}_M^f$ as constant linear
combinations of $\partial_a=\partial /\partial x^a.$
As a notation for a flat structure, we will use indiscriminately
$\nabla_0$ or $\Cal{T}_M^f$. 

\smallskip

Now consider a manifold $M$ whose structure
sheaf is endowed with an $\Cal{O}_M$--bilinear
(super)commutative
and associative multiplication $\circ$,
and eventually with identity $e$. 
In the following definition we do not assume that it satisfies
(1.2).

\medskip

{\bf 2.2. Definition.} {\it a) A flat structure $\Cal{T}_M^f$
on $M$ is called compatible with $\circ$, 
if in a neighborhood of any point there exists a vector field
$C$ such that for arbitrary local flat vector fields $X,Y$
we have
$$
X\circ Y= [X,[Y,C]].
\eqno(2.1)
$$
$C$ is called a local vector potential for $\circ$.

\smallskip
b)  $\Cal{T}_M^f$ is called compatible with $(\circ ,e)$, if
a) holds and moreover,  $e$ is flat.}

\smallskip

{\bf 2.3. Remarks.} (i) If we choose a
local flat coordinate system $(x^a)$ and write 
$C=\sum_c C^c\partial_c$, $X=\partial_a,\,Y=\partial_b$, then (1.1) becomes 
$$
\partial_a\circ\partial_b =\sum_c C_{ab}{}^c\partial_c,\quad
C_{ab}{}^c= \partial_a\partial_bC^c . 
\eqno(2.2)
$$
If moreover
$e$ is flat, we may choose local flat
coordinates so that  $e=\partial_0$, and the
conditions $e\circ \partial_b=\partial_b$ reduce to 
$C_{0b}{}^c=\delta_b^c.$

\smallskip

(ii) If we choose an arbitrary $C$ and {\it define} a composition
$\circ :\,\Cal{T}_M^f\otimes \Cal{T}_M^f\to \Cal{T}_M$ by the formula
(2.1), it will be automatically supercommutative
in view of the Jacobi formula. Associativity, however, is a quadratic differential constraint on $C$ which was called ``oriented associativity
equations'' in [LoMa2].

\smallskip

(iii) If $C$ exists, it is not unique. As one sees from (2.2), locally
it is defined modulo the span of vector fields
$\{\partial_a, x^b\partial_c\}$ which form a subsheaf of Lie algebras in 
$\Cal{T}_M$ depending on $\nabla_0$ and denoted $\Cal{T}_M^{(1)}$. 

\medskip

{\bf 2.4. Proposition.} {\it Assume that a multiplication $\circ$
on the tangent sheaf of $M$ admits a compatible flat structure.
Then it satisfies (1.2). Thus, $(M,\circ )$ is an $F$--manifold.}

\medskip

{\bf Proof.} The calculation is essentially the same as in
the proof of Theorem 2 of [HeMa]. We reproduce it for
completeness and because the context is slightly different.

\smallskip

First of all, the identity (1.2) can be rewritten
as follows: for any local vector fields
$X,Y,Z,W$ we have
$$
[X\circ Y,Z\circ W]-[X\circ Y,Z]\circ W-
(-1)^{(X+Y)Z}Z\circ [X\circ Y,W]
$$
$$
-X\circ [Y,Z\circ W] -(-1)^{XY}Y\circ [X,Z\circ W]+ X\circ [Y,Z]\circ W
+(-1)^{YZ} X\circ Z\circ [Y,W]
$$
$$
 + (-1)^{XY}Y\circ [X,Z]\circ W
+(-1)^{X(Y+Z)}Y\circ Z\circ [X,W] =0
\eqno(2.3)
$$
This form is convenient because it turns out that the left hand
side of (2.3) is in fact a tensor, that is $\Cal{O}_M$--polylinear
in $X,Y,Z,W$: cf. [Me1], [Me2] for a discussion and a generalization of this 
identity.

\smallskip

Therefore to verify (2.3) in our context it suffices
to check that the left hand side vanishes on all quadruples
of local flat fields $(\partial_a, \partial_b, \partial_c, \partial_d).$
Since flat fields (super)commute, 
the last four summands of (2.3) vanish, and only the first five ones
should be taken care of. Let us denote the structure ``constants'' $C_{ab}{}^c$ as in (2.2). Calculating the coefficient of $\partial_f$ in 
the left hand side of (2.3),
we represent it as a sum of the respective five summands,
for which we introduce a special notation in order
to explain the pattern of cancellation:
$$
\sum_e C_{ab}{}^e\partial_eC_{cd}{}^f-(-1)^{(a+b)(c+d)}
\sum_e C_{cd}{}^e\partial_eC_{ab}{}^f=\alpha_1+\beta_1,
$$
$$
(-1)^{(a+b)c}\sum_e \partial_cC_{ab}{}^eC_{ed}{}^f=\alpha_2,\quad
(-1)^{(a+b+c)d}\sum_e \partial_dC_{ab}{}^ecC_{ec}{}^f=\gamma_1,
$$
$$
-(-1)^{a(b+c+d)}\sum_e \partial_bC_{cd}{}^eC_{ea}{}^f=\gamma_2,\quad
-(-1)^{(c+d)b}\sum_e \partial_aC_{cd}{}^eC_{eb}{}^f=\beta_2.
$$
Here we write, say, $(-1)^{(a+b)c}$ as a shorthand for
$(-1)^{(\widetilde{x_a}+\widetilde{x_b})\widetilde{x_c}}$.

\smallskip

Using the second formula (2.2), we can replace 
$\partial_eC_{cd}{}^f$ in $\alpha_1$ by
$(-1)^{ec}\partial_cC_{ed}{}^f$.
After this we see that
$$
\alpha_1+\alpha_2=(-1)^{(a+b)c}\partial_c
\left(\sum_e C_{ab}{}^eC_{ed}{}^f\right).
$$
Similarly, permuting $a$ and $e$ in $\beta_1$ we find
$$
\beta_1+\beta_2=-(-1)^{(c+d)b}\partial_a
\left(\sum_e C_{cd}{}^eC_{eb}{}^f\right).
$$ 
Now rewrite $\gamma_1$ permuting $a,d,$ and $\gamma_2$
permuting $b,c.$ Calculating finally 
$\beta_1+\beta_2+\gamma_1+\gamma_2$ we see that it cancels with
$\alpha_1+\alpha_2$ due to the associativity of $\circ$ written as 
in [Ma], I (1.5).

\bigskip

{\bf 2.5. Pencils of flat connections.} Let $M$ be a supermanifold
endowed with a torsionless flat connection $\nabla_0:\,\Cal{T}_M\to
\Omega^1_M\otimes_{\Cal{O}_M}\Cal{T}_M$ and an odd global section 
$A\in \Omega^1_M\otimes_{\Cal{O}_M} End(\Cal{T}_M)$ (it is called a Higgs field
in other contexts).

\smallskip 

We will use this operator in order to define
two structures. 

\smallskip

First, for any even constant $\lambda$,
consider a connection $\nabla_{\lambda}=\nabla_{\lambda}^A$ on $\Cal{T}_M$:
$$
\nabla_{\lambda}:= \nabla_0 + \lambda \,A.
\eqno(2.4)
$$
Write the curvature form of $\nabla_{\lambda}$ as
$$
\nabla_{\lambda}^2=\lambda R_1+\lambda^2R_2.
$$
Second, define an $\Cal{O}_M$--bilinear composition law 
$\circ =\circ^A$ on $\Cal{T}_M$:
$$
X\circ^A Y:=i_X(A)(Y),\quad i_X(df\otimes G):=Xf\cdot G.
\eqno(2.5)
$$ 

\medskip

{\bf 2.6. Proposition.} {\it a) $\circ^A$ is supercommutative if and only
if all connections $\nabla_{\lambda}^A$ are torsionless.

\smallskip

b) Assume that a) holds. Then $\circ^A$ restricted  
to $\nabla_0$--flat vector fields
can be written everywhere locally in the form (2.1)
if and only if $R_1=0$.

\smallskip

c) Assume that a) holds.  Then $\circ^A$ is associative
if and only if $R_2=0$.}

\medskip

{\bf 2.7. Corollary.} {\it $(M,\circ^A,\nabla_0)$
is an $F$--manifold with compatible flat structure if and
only if $\nabla_{\lambda}^A$ is a pencil
of torsionless flat connections.

\smallskip

In this case, $(M,\circ^A,\nabla_{\lambda}^A)$
is an $F$--manifold with compatible flat structure
for any $\lambda$ as well.}

\medskip

{\bf Proof.} a) Omitting for brevity the superscripts $A$,
we can write covariant derivatives as
$$
\nabla_{\lambda,X}(Y) =\nabla_{0,X}(Y)+\lambda X\circ Y.
$$

\smallskip

We will again work with a local basis of $\nabla_0$--flat
vector fields $\{\partial_a\}$. Then the vanishing of torsion of $\nabla_{\lambda}$ means that
$$
\nabla_{\lambda ,\partial_a}(\partial_b)=(-1)^{ab}
\nabla_{\lambda ,\partial_b}(\partial_a)
$$
that is
$$
\lambda \,\partial_a\circ\partial_b = (-1)^{ab}\lambda\,
\partial_b\circ\partial_a
$$
which is the supercommutativity of $\circ$.

\smallskip

b) Define the structure/connection
coefficients $A_{ab}{}^c$ by the formula 
$\partial_a\circ\partial_b =\sum_c A_{ab}{}^c\partial_c.$ 
Since we have assumed a), these coefficients
are symmetric in $a,b$.
The vanishing of $R_1$ is equivalent to
$$
\forall\, a,b,c,e,\quad \partial_a A_{bc}{}^e=(-1)^{ab}\partial_b A_{ac}{}^e.
$$
This condition means that local 1--forms
$\sum_b dx^bA_{bc}{}^e$ are closed. In view of the De Rham lemma,
they are locally exact,
which is equivalent to the existence 
of local functions $B_c{}^e$ such that
$$
A_{bc}{}^e=\partial_b B_c{}^e.
$$
Again, since  $A_{bc}{}^e$ is symmetric in $b,c$,
the form $\sum_cdx^cB_c{}^e$ is closed, and the De Rham lemma
implies existence of some local
functions $C^e$ such that $B_c{}^e=\partial_cC^e$.
Hence the local vector field $C:=\sum C^e\partial_e$
determines our $\circ$ as in formula (2.1). Retracing these
arguments in reverse order, one sees that if $\circ$ is locally of the form
(2.1), then $R_1=0$. 

\smallskip

c) Commutativity assumed, the equivalence of 
associativity with vanishing of $R_2$
is checked exactly as in the proof of Theorem 1.5. b) in [Ma], p. 21,
and we will not repeat it. This completes the proof
of the Proposition 2.6.

\smallskip

The first statement of Corollary 2.7  is thereby proved as well.

\smallskip
The second statement follows from the fact that
changing the initial point on the affine line of flat
torsionless connections does not affect the vanishing
of torsion and curvature.

\medskip

{\bf 2.8. Auxiliary formulas.} Let $(M,\circ ,\nabla )$
be an $F$--manifold with a compatible flat structure.
Here and below we write $\nabla$ in place of former $\nabla_0$.
We will often use formulas $\nabla_XY-\nabla_YX=[X,Y]$
(vanishing torsion) and $\nabla_{[X,Y]}=[\nabla_X,\nabla_Y]$
(vanishing curvature). If $X$ is flat, $\nabla_XY=[X,Y].$
Put 
$$
D(X,Y,Z):=\nabla_X(Y\circ Z)-\nabla_X(Y)\circ Z-(-1)^{XY}Y\circ\nabla_XZ
\eqno(2.6)
$$
(compare this with $P_X(Y,Z)$ defined by (1.1) and depending
only on $\circ$).

\medskip

{\bf 2.8.1. Lemma.} {\it $D(X,Y,Z)$ is a symmetric tensor
with values in $\Cal{T}_M$.}

\smallskip

{\bf Proof.} One first directly checks that 
$D(X,Y,Z)-(-1)^{XY}D(Y,X,Z)$ is a tensor. Then to show that 
this difference vanishes 
identically it suffices to prove this for flat $X,Y,Z$.
Only two terms of six survive, and they can be rewritten 
in view of (2.1) as
$$
\nabla_X [Y,[Z,C]]-(-1)^{XY}\nabla_Y [X,[Z,C]]=
[X, [Y,[Z,C]]]-(-1)^{XY}[Y,[X,[Z,C]]]
$$
which vanish thanks to Jacobi. Moreover, $D(X,Y,Z)$ is obviously
symmetric in $Y,Z$, so symmetric in all three
arguments. Finally, it is $\Cal{O}_M$--linear in $X$, which
completes the argument.

\medskip

{\bf 2.9. Nonflat identities.} In this subsection, we collect for further use
several formulas involving an $F$--manifold $(M,\circ ,\nabla ,e)$ 
with a compatible flat structure but nonnecessarily flat
identity. 

\smallskip

Denote by $\Cal{L}$ the set of global (or
local) vector fields $\varepsilon$ satisfying the following
condition:

\smallskip 

{\it (*) for any local vector field $Y$ we have}
$$
\nabla_Y\varepsilon =Y\circ \nabla_e\varepsilon .
\eqno(2.7)
$$

Clearly, $\Cal{L}$ is a (sheaf of) vector space(s). 
Notice that we could have postulated a formally weaker condition:
$\nabla_Y\varepsilon =Y\circ l(\varepsilon )$ for some
functional $l$, but putting $Y=e$ we immediately get
the unique possibility $l(\varepsilon )=\nabla_e\varepsilon$.  

\smallskip

Another useful remark is this: if (2.7) holds for $\nabla$
($\varepsilon$ being fixed), then
it holds for any $\nabla_{\lambda}$ in the relevant
pencil of flat connections:
$$
\nabla_{\lambda ,Y}\varepsilon = \nabla_Y\varepsilon + \lambda Y\circ\varepsilon
=Y\circ \nabla_{\lambda ,e}\varepsilon .
$$
More generally, if (2.7) holds for any single $\nabla_{\lambda}$,
it holds for all of them.
 
\smallskip

{\bf 2.9.1. Proposition.} {\it a) $\Cal{L}$ contains $e$,  $\roman{Ker}\,\nabla_{\lambda}$ for all $\lambda$,
and is closed with respect to $\nabla_e$ so that it contains
the span of vector fields $\{\roman{Ker}\,\nabla_{\lambda}, e, \nabla_e{e}, \nabla_e^2{e},
\dots \}$.

b) For any $\varepsilon \in \Cal{L}$ we have
$$
\roman{ad}\,\varepsilon = \nabla_{\varepsilon} - \circ\nabla_e
\varepsilon
\eqno(2.8)
$$
and
$$
P_{\varepsilon}(Y,Z)=D(\varepsilon ,Y,Z)+Y\circ Z
\circ\varepsilon.
\eqno(2.9)
$$

c) The operator $\roman{ad}\,e$ is a $\circ$--derivation:
$$
[e,Y\circ Z]=[e,Y]\circ Z+Y\circ [e,Z].
\eqno(2.10)
$$
}

{\bf Proof.} Clearly, $\Cal{L}$ contains $\roman{Ker}\,\nabla_{\lambda}$. We have from (2.6)
$$
D(e,e,Y)=-\nabla_ee\circ Y,\quad D(Y,e,e)=-\nabla_Ye,
$$
These fields coincide since $D$ is symmetric, which shows
that $e\in \Cal{L}$. 

Since generally
$[\varepsilon ,Y]=\nabla_\varepsilon Y-\nabla_Y\varepsilon$, 
(2.8) is a consequence of (2.7). Inserting (2.8) in the definition
of $P_{\varepsilon}$, we get (2.9).
 
Putting $X=Y=e$ in (1.2), we obtain (2.10). 

Now we will check that $\Cal{L}$ is $\nabla_e$--stable. In fact,
let $\varepsilon \in \Cal{L}$. Then
$$
\nabla_Y\nabla_e\varepsilon=
\nabla_e\nabla_Y\varepsilon +\nabla_{[Y,e]}\varepsilon=
\nabla_e (Y\circ \nabla_e\varepsilon )+ [Y,e]\circ\nabla_e\varepsilon .
$$
In the last expression, replace $\nabla_e$ by $\roman{ad}\,e+\circ\nabla_ee$
and use the fact that $\roman{ad}\,e$ is a $\circ$--derivation. We get
$$
[e,Y\circ\nabla_e\varepsilon ]+Y\circ \nabla_e\varepsilon\circ\nabla_ee+
[Y,e]\circ\nabla_e\varepsilon =
Y\circ [e,\nabla_e\varepsilon ] + Y\circ \nabla_e\varepsilon\circ\nabla_ee
$$
which has the required form (2.7) so that $\nabla_e\varepsilon\in\Cal{L}$.

\medskip

{\bf 2.10. Identities with $\nabla_ee=ce$.} Start with an $F$--manifold
with a compatible flat structure and flat identity
$(M,\circ,e,\nabla_0)$. Replace $\nabla_0$ by some shifted
connection $\nabla :=\nabla_0+cA$
in the pencil of flat connections (2.4).
In view of Corollary
2.7, it will still be compatible with $\circ$, but
$e$ will be nonflat: $\nabla_ee=c e.$ 

\smallskip

Conversely, if $\nabla_ee=ce$ for some constant $c$, then we can choose
a new  initial point $\nabla_0$ in the affine line of flat torsionless connections
$\nabla +\lambda A$ compatible with $\circ$ in such a way that $e$ will become $\nabla_0$--flat:
simply put $\lambda = -c.$

\smallskip

Examples below show that an $F$--manifold may admit a 
$\circ$--invariant flat metric
whose Levi--Civita connection belongs to the
same affine line of flat torsionless connections compatible with $\circ$,
but does not coincide with the point of flat identity.

\medskip

{\bf 2.10.1. Examples.} (i) Here we will demonstrate that Dubrovin's 
``almost Frobenius structure''
([Du2], Def. 9) after forgetting the metric but retaining the affine 
flat structure $\nabla$ corresponding to its Levi--Civita connection,
 becomes an $F$--manifold
with compatible flat structure and nonflat identity $e$
satisfying the condition above: $\nabla_ee=ce$ with a constant $c$.

\smallskip

In fact, the Levi-Civita connection $\nabla$ is torsionless and flat.
Compatibility condition (2.1) follows from Dubrovin's
potentiality formulas (3.20) and (3.21) of [Du2]. Finally,
in a system of flat coordinates $(x^a)$ (Dubrovin's $(p^i)$),
the identity $e$ takes form $\dfrac{1-d}{2} \sum_ax^a\partial_a$
(Dubrovin's identity/Euler field $E$, formula (3.23) of [Du2]), 
where $d\ne 1$ is a constant.

\smallskip

Hence $\nabla_ee=ce$ where $c=\dfrac{(1-d)^2}{4}.$

\smallskip 

(ii) Quantum $K$--theory formal $F$--manifolds studied 
in [Gi] and [Lee] share a similar property: $\nabla_e e = e/2$.
I am thankful to Y.-P. Lee who has shown me a  calculation
establishing this (it was known to Givental). 
It follows directly from the form of Christoffel
symbols, presented in [Lee], sec. 5.2.
 
\medskip

{\bf 2.11. Summary.} The ``space'' of flat affine connections compatible
with a given multiplication $\circ$ is fibered by affine lines.
Without an additional structure, there is generally no way
to single out a point on such a line. Point of a flat unity
and point of the Levi--Civita connection of an invariant flat
metric may in general diverge, although they
coincide for Frobenius manifolds in the strict sense. 

\smallskip

Therefore it makes
sense to summarize the resultats of this
section stressing  the $\lambda$--coordinate free
aspect. 

\smallskip

Let us start with a structure $(M,\Cal{P})$ where
$M$ is a (super)manifold, and $\Cal{P}$ is an affine line (pencil)
of  connections on $\Cal{T}_M$.
Such a line determines the following derivative structures on $M$:

\medskip

(i) {\it A $\bold{C}^*$--torsor $\Cal{A}$ of odd tensors 
$A\in \Omega^1_M\otimes_{\Cal{O}_M} End\,(\Cal{T}_M)$.}

\smallskip

Namely, each $A\in \Cal{A}$ is a difference $\nabla_1-\nabla_2$
of some $\nabla_1\ne\nabla_2\in \Cal{P}$, and for any two such
differences $A,\,B$ we have $A=\lambda B$ for some $\lambda
\in \bold{C}^*$.

\medskip

(ii)  {\it A $\bold{C}^*$--torsor $\Cal{M}$ of 
 $\Cal{O}_M$--bilinear multiplications $\circ^A$ on $\Cal{T}_M$:
for each $A\in \Cal{A}$, we define $X\circ^A Y:=i_X(A)(Y).$}

\smallskip

One (and hence all) $\circ^A$ is (are) supercommutative,
if and only if  two (and hence all) $\nabla\in \Cal{P}$ are torsionless.
Assume that this condition is satisfied. 

\smallskip

One (and hence all) $\circ^A$ is (are) associative
and admit local vector potentials, 
if and only if  all $\nabla\in \Cal{P}$ are flat.

\smallskip

If one of the multiplications $\circ^A$ has an identity $e_A$,
then each of them has an identity: $A=\lambda B$ implies
$X\circ^AY=\lambda X\circ^BY$, so that $e_A=\lambda^{-1}e_B$.
Hence this is a property of the whole pencil $\Cal{P}$,
which we will then call {\it unital}. All identities $e_A$
form a $\bold{C}^*$--torsor as well.

\smallskip

For an unital pencil $\Cal{P}$ there may exist
at most one $\nabla_0 \in \Cal{P}$ such that 
some (equivalently, any) $e_A$ is $\nabla_0$--flat.
It exists if and only if $e_A$ belongs
to an eigenspace of some (equivalently, any)
$\nabla_e$,  $\nabla \in \Cal{P}$.
Such $\Cal{P}$ can be called {\it flat unital.}

\medskip

(iii) {\it The sheaf $\Cal{L}$ of such local vector fields
$\varepsilon$ that for any local vector field $Y$,
some (equivalently, any) $\nabla\in \Cal{P}$,
and some (equivalently, any) $A\in \Cal{A}$ there exists
a vector field $l(\varepsilon )$ such that
$\nabla_Y\varepsilon =Y\circ l(\varepsilon )$.}

\smallskip

From Proposition 2.9.1 we know that $\Cal{L}$
contains $\sum_{\nabla\in\Cal{P}} Ker\,\nabla$.
If $\Cal{P}$ is unital, then $\Cal{L}$ contains all 
identities $e$ as well and is stable with respect to $\nabla_e$
for any $\nabla \in \Cal{P}.$

\bigskip

\centerline{\bf \S 3. Extended flat connections and Euler fields}

\medskip

{\bf 3.1. Notation.} We continue considering an $F$--manifold
with a compatible flat structure and possibly nonflat identity
$(M,\circ,\nabla ,e)$. Let $\lambda$
be the coordinate on the affine line of flat
metrics as in (2.4). We will generally denote by
$\widehat{M}$ the total space of the constant family of manifolds $M$ 
over a base which
might be the $\lambda$--affine line, or its projectivization,
or formal completion, say, at $\lambda =\infty$ etc,
to be specified in concrete situations. Let $pr_M:\, \widehat{M}\to M$
be the projection on the base.

\smallskip

We will identify $\Cal{T}_M$ with the subsheaf 
of $\lambda$--independent vector fields in $\Cal{T}_{\widehat{M}}$ .
We denote by $\partial_{\lambda}$ the vector field on
$\widehat{M}$ which annihilates $\Cal{O}_M$ and
such that $\partial_{\lambda}\lambda =1.$  

\smallskip

We will now consider $\widehat{\nabla}$ defined by covariant derivatives
along $\lambda$--indepen\-dent vector fields 
$$
\widehat{\nabla}_X:= \nabla_X +\lambda X\circ
\eqno(3.1)
$$ 
(cf. (2.4)) as a part of a
connection on $pr^*_M(\Cal{T}_M)$ over $\widehat{M}$
(it might be meromorphic or formal).
To complete (3.1) to a full connection $\widehat{\nabla}$
on $pr^*_M(\Cal{T}_M)$ ({\it not} on $\Cal{T}_{\widehat{M}}$!), we should
choose a covariant deriva\-tive along $\partial_{\lambda}$
which on $Y\in pr^*_M(\Cal{T}_M)$ we require to be of the form
$$
\widehat{\nabla}_{\partial_{\lambda}}(Y) = \partial_{\lambda}Y+H(Y)
\eqno(3.2)
$$
where $H$
is an even endomorphism of $pr_M^{*}(\Cal{T}_M))$ 
which we will have to allow
to depend rationally or even  formally of $\lambda$.
We will write $\widehat{\nabla}^H$ for $\widehat{\nabla}$
if we need to stress its dependence on $H$.

\medskip

{\bf 3.2. Flatness conditions.} We want $\widehat{\nabla}^H$
to be flat. In this subsection we will spell
the implied conditions on $H$, which are equivalent to
$[\widehat{\nabla}_{X},\widehat{\nabla}_{\partial_{\lambda}}]=0$
for all $\lambda$--independent $X$. It suffices to check
this identity by applying it to $\lambda$--independent fields $Y$.
We have
$$
\widehat{\nabla}_{X}\widehat{\nabla}_{\partial_{\lambda}} (Y)=
\nabla_{X} H(Y)+\lambda X\circ H(Y)
\eqno(3.3)
$$
where $\nabla_X$ and $\circ$ are extended to $pr_M^{*}(\Cal{T}_M)$
in the evident way. Furthermore,
$$
\widehat{\nabla}_{\partial_{\lambda}}\widehat{\nabla}_{X}(Y)=
\widehat{\nabla}_{\partial_{\lambda}}(\nabla_XY+\lambda X\circ Y)=
$$
$$
X\circ Y +H(\nabla_XY) +\lambda H(X\circ Y).
\eqno(3.4)
$$
Combining (3.3) and (3.4), we get the following reformulation of the flatness condition:
$$
\forall\, X,Y,\quad 
 H(X\circ Y) =X\circ H(Y)+\frac{1}{\lambda}\left(\nabla_XH(Y)
-X\circ Y-H(\nabla_XY)\right).
\eqno(3.5)
$$
\smallskip
Putting here $Y=e$ and using $\nabla_Xe=X\circ\nabla_ee$
we get a functional equation for $H$
which is in principle {\it  only necessary} for $\widehat{\nabla}$
to be flat:
$$
H(X)=X\circ H(e)+\frac{1}{\lambda}\left(\nabla_XH(e)
-X-H(X\circ\nabla_ee)\right).
\eqno(3.6)
$$
If (3.6) is satisfied, we can put here $X=e$ and get
$\nabla_eH(e)=e+H(\nabla_ee)$, or better, since 
$\nabla_eH(e)=[e,H(e)]+\nabla_{H(e)}e=[e,H(e)]+H(e)\circ\nabla_ee$
(see (2.7) and Prop. 2.12.1),
$$
[e,H(e)]+H(e)\circ\nabla_ee-H(\nabla_ee)=e.
\eqno(3.7)
$$
\smallskip

{\bf 3.3. Formal solutions to (3.6).} Let $T$ be the space of
global vector fields on $M$. Put $\mu := \lambda^{-1}$,
$e_1:=\nabla_ee$, and 
$$
(e+\mu e_1)^{-1}:=\sum_{i=0}^{\infty}
(-1)^i(e_1)^{\circ i}\mu^i \in T[[\mu ]].
\eqno(3.8)
$$

\smallskip

{\bf 3.3.1. Theorem.} {\it a) There exists an 1--1 correspondence
between the solutions to (3.6) $H\in \roman{End}\,T[[\mu ]]$
and the solutions $E\in T[[\mu]]$ to the equation
$$
(e+\mu e_1)\circ\nabla_{(e+\mu e_1)^{-1}}E-  e_1\circ E = e.
\eqno(3.9)
$$
Namely, from $E$ one reconstructs $H$ as follows:
$$
H(X):=X\circ E+\mu (\nabla_{(e+\mu e_1)^{-1}\circ X}E- (e+\mu e_1)^{-1}
\circ X)+
((e+\mu e_1)^{-1}\circ X-X)\circ E ,
\eqno(3.10)
$$
and from $H$ one gets $E$ as $E=H(e)$.

\bigskip

b) A solution $H=H_E$ as in a) satisfies the flatness condition (3.5)
if and only if for all $\nabla$--flat vector fields $X,Y$ we have
$$
P_E(X,Y)-[X,[Y,C-\mu E]] = X\circ H(Y\circ e_1)+\mu [X,H(Y\circ e_1)],
\eqno(3.11)
$$
where $C$ is defined by $X\circ Y=[X,[Y,C]]$ as in (2.1).
}
\bigskip

{\bf Proof.} a) Assume first that $H$ satisfies (3.6). Put
$E:=H(e)$. We will prove that (3.10) holds.
In fact, (3.6) can be rewritten as
$$
H(X)=A(X)+\mu B(X) +\mu H(C(X))
\eqno(3.12)
$$
where $A,B,C$ are $\Cal{O}_M[[\mu]]$--linear operators
$$
A(X)=X\circ E,\quad B(X)=\nabla_XE-X,\quad C(X)=-X\circ e_1.
\eqno(3.13)
$$
We can replace in (3.12) $X$ by $C(X)$ and then put the
obtained formula for $H(C(X))$ into the right hand side
of (3.12). Infinitely iterating this procedure, we obtain:
$$
H(X)=A(X)+\sum_{n=1}^{\infty}\mu^nBC^{n-1}(X)+\sum_{n=1}^{\infty}\mu^nAC^{n}(X)=
$$
$$
X\circ E+\mu B((e+\mu e_1)^{-1}\circ X)+A((-e+(e+\mu e_1)^{-1})\circ X)
$$
which is equivalent to (3.10).

\smallskip

Now put $X=e$ in (3.10), and multiply the result by $\mu^{-1}(e+\mu e_1)$.
We will get (3.9).

\smallskip

Conversely, start with an arbitrary $E\in T[[\mu ]]$ and define 
$H(X)$ by the formula (3.10). Retracing backwards the calculations above,
one sees that $H(X)$ satisfies (3.12). This means that a version of (3.6) holds
in which $H(e)$ is replaced by $E$. The additional condition $E=H(e)$ is then
equivalent to (3.9). This completes the proof.

\smallskip

b) Now consider the full flatness condition (3.5).
Both sides of 
(3.5) are $\Cal{O}_M$--bilinear in $X,Y$ so that it suffices
to check its meaning for $\nabla$--flat $X,Y$, and we will assume this
in the following calculations. The term $H(\nabla_XY)$ in (3.5) will vanish.
We replace $\nabla_X$ by $[X,*]$ and similarly for $\nabla_Y$.

\smallskip
 
Replace  $X$ in (3.6) first by $X\circ Y$ and then by $Y$.
Put the resulting expressions for $H(X\circ Y)$ and $H(Y)$ 
into (3.5). After some cancelations and division by $\mu$ we find
$$
\forall\ \roman{flat}\ X,Y,\quad 
\nabla_{X\circ Y}E=
$$
$$
X\circ [Y,E]+[X,Y\circ E] -X\circ Y+
\mu [X,[Y,E]]-X\circ H(Y\circ e_1)-\mu [X,H(Y\circ e_1)]
\eqno(3.14)
$$
Now rewrite the left hand side of (3.14) using (2.6) and Lemma 2.11.1: 
$$
\forall\ \roman{flat}\ X,Y,\quad
\nabla_{X\circ Y}E=[X\circ Y,E]+\nabla_E(X\circ Y)=
[X\circ Y,E] +D(E,X,Y)=
$$
$$
[X\circ Y,E] +D(X,Y,E)=
[X\circ Y,E]+[X,Y\circ E]-(-1)^{XY}Y\circ [X,E].
$$
Putting this into (3.14), after cancelations and regrouping we obtain (3.11).

\medskip

{\bf 3.4. Extended connections and Euler fields.} In this subsection we 
assume that
$e_1=ce$ for a constant $c$ (as in 2.13 above) and moreover, that $E:=H(e)$ 
is independent of $\lambda$. Shifting the base point in the relevant
pencil of flat connections we may and will assume that $c=0$.
Equivalently, we should take $\mu /(1+c\mu)$ for new $\mu$,
and of course, adjust the notion of flatness.
In this case (3.6) and (3.7) reduce respectively to
$$
H(X)=X\circ E+ \mu \left(\nabla_XE-X\right),
\eqno(3.15)
$$
$$
[e,E]=e.
\eqno(3.16)
$$ 

We now recall the definition from [HeMa].

\smallskip

{\bf 3.4.1. Definition.} {\it a) A  vector field $E$ on an
$F$--manifold $(M,\circ )$ is called an Euler field of
weight $d_0$ if
$$
P_E(X,Y)=d_0 X\circ Y.
\eqno(3.17)
$$ 
b) The compatibility condition of $E$ with a compatible flat structure
$\Cal{\nabla}$ reads
$$
[E,\roman{Ker}\,\nabla ] \subset \roman{Ker}\,\nabla .
\eqno(3.18)
$$}

Local Euler fields on $(M,\circ )$ form a sheaf of 
linear spaces and Lie algebras.
Weight is a linear function on this sheaf. Commutator
of two Euler fields is an Euler field of weight zero.
Identity is an Euler field of weight zero.

\smallskip

These statements hold for Euler fields
compatible with a given $\nabla$ as well. Identity
as an Euler field is compatible with $\nabla$ 
if it is $\nabla$--flat or more generally
$e\in \Cal{T}_M^f$.

\smallskip

{\bf 3.4.2. Proposition.} {\it The connection $\widehat{\nabla}^H$
with $H$ of the form (3.15) and $\lambda$--independent $E$
is flat if and only if $E$ is an Euler field of weight 1
compatible with $\nabla$.}

\medskip

{\bf Proof.} In the general flatness condition (3.11) the right hand
side vanishes. Since $E$ is $\lambda$--independent,
the term $[X,[Y,\mu E]]$ in the left hand side must vanish
as well which means that $[E,\Cal{T}^f_M]\subset \Cal{T}^f_M$.
After that (3.11) reduces to (3.17) with $d_0=1$.
This completes the proof.

\medskip

{\bf 3.4.3. A further deformation of $\widehat{\nabla}$.} 
Assume that we are given 
 an Euler field $E$ of weight 1
compatible with $\nabla$ and that $e\in \Cal{T}_M^f$.
Then we can construct a line of Euler fields
$E^{(s)}:= E+se$ of weight 1
compatible with $\nabla$, and the respective deformation
$\widehat{\nabla}^{(s)}$ of $\widehat{\nabla}$
obtained by using $E^{(s)}$ (or rather $E^{(s\mu)}$)
in place of $E$ in (3.15). 

\smallskip

See a discussion of this family in the context
of Frobenius manifolds on pp. 154--157 of [He].

\bigskip

\centerline{\bf \S 4. Dubrovin's duality}

\medskip

{\bf 4.1. Pencils of flat connections on an external bundle.}
Consider now a slight variation of the setup
described in 2.11: a locally free sheaf $\Cal{F}$ on a manifold
$M$ and a pencil (affine line) $\Cal{P}$ of connections
$\nabla:\,\Cal{F}\to
\Omega^1_M\otimes_{\Cal{O}_M}\Cal{F}$. The difference
of any two connections is now an odd global section
$\nabla_1-\nabla_0=A\in \Omega^1_M\otimes_{\Cal{O}_M} End\,\Cal{F}$.
Any two differences are proportional; we will often choose
one arbitrarily.
Any $\nabla$ can be extended to an odd derivation 
(denoted again $\nabla$) of $\Omega^*_M\otimes_{\Cal{O_M}}T(\Cal{F})$
in the standard way, where $T(\Cal{F})$ is the total
tensor algebra of $\Cal{F}$.

\smallskip

From now on, we will assume that all connections in $\Cal{P}$ are flat.
This means that for any $\nabla \in \Cal{P}$, $\nabla A=0$ and $A\wedge A =0$.
Again, in view of the De Rham lemma, $A$ can be everywhere
locally written as $\nabla B$ where $B=B_{\nabla}$ (it generally depends
on $\nabla$) is an even section
of $End\,\Cal{F}$. Notice that we cannot assume $\Cal{P}$
to be torsionless: this makes no sense for an external bundle.
We will see that this restriction can be replaced by
fixing an additional piece of data.

\smallskip

Let now $U$ be a coordinate neighborhood in $M$ over which
the linear superspace $F$ of local $\nabla$--flat sections of
$\Cal{F}$ trivializes $\Cal{F}.$ Denote by $\widetilde{F}$
the affine supermanifold associated with $F$. 
Let $q:\,\widetilde{\Cal{F}}\to M$ be the fibration
of supermanifolds which is ``the total space'' of
$\Cal{F}$ as a vector bundle: for example, in
algebraic geometry this is the relative affine spectrum
of $\roman{Symm}_{\Cal{O}_M}\,(\Cal{F}^*)$, $\Cal{F}^*$
being the dual sheaf. Hence local sections of $\Cal{F}$
become local sections of $q$.

\smallskip

Clearly, $\nabla$ trivializes $q$ over $U$: we have
a well defined isomorphism $q^{-1}(U)=\widetilde{F}\times U$
turning $q$ into projection. 
 
\smallskip

Let now $u\in F$ be a $\nabla$--flat section of $\Cal{F}$
over $U$, $\nabla B=A$, $B\in \roman{End}\,\Cal{F}.$ 
Then $Bu$ is a section of
$\Cal{F}$; we will identify it with a section of $q$
as above. Projecting to $\widetilde{F}$,
we finally get a morphism $Bu:\,U\to \widetilde{F}$.
Thus $Bu$ denotes several different although closely related objects;
hopefully, this will not lead to a confusion.

\smallskip

In a less fancy language, if we choose a basis of flat sections in $F$
and a system of local coordinates $(x^a)$ on $M$,
$B$ becomes a local even matrix function $B(x)$ acting from the left
on columns of local functions. Then $Bu$ becomes
the map $U\to \widetilde{F}:\, x\mapsto B(x)u$.
Since $B$ is defined up to $\roman{Ker}\,\nabla$,
this map is defined up to a constant shift.

\medskip

{\bf 4.2. Definition.} {\it A section
$u$ of $\Cal{F}$ is called a primitive section
with respect to  $\nabla \in\Cal{P}$ ,
if it is $\nabla$--flat, and $Bu$ is a local
isomorphism of $U$ with a subdomain of
$\widetilde{F}$.}

\smallskip

Since $F$ is linear, $\widetilde{F}$ has a canonical flat structure
$\nabla^F:\,\Cal{T}_{\widetilde{F}}\to \Omega^1_{\widetilde{F}}
\otimes \Cal{T}_{\widetilde{F}}.$ If $u$ is primitive,
$Bu$ is a local isomorphism allowing us to identify locally
$(Bu)^*(\Cal{T}_{\widetilde{F}}) = \Cal{T}_M$. Moreover, we can
consider the pullback of $\nabla^F$ with respect to $Bu$:
$$
\nabla^*:= (Bu)^*(\nabla^F): \Cal{T}_M\to \Omega^1_M\otimes\Cal{T}_M.
\eqno(4.1)
$$ 
From the remarks above it is clear that the local flat structures
induced by the maps $Bu$ on $M$ do not
depend on local choices of $B$ and glue together to
a flat structure on all of $M$ determined by $\nabla$ and $u$.
 
\smallskip

It follows also that $\nabla^*$ is flat and torsionless.

\smallskip

There is another important isomorphism produced by this
embedding. Namely, restricting $(Bu)^*$ to $F\subset \Cal{T}_{\widetilde{F}}$
we can identify $(Bu)^*(F)\subset \Cal{T}_M$ with $F\subset \Cal{F}$
and then by $\Cal{O}_M$--linearity construct an
isomorphism $\beta^*:\Cal{F}\to\Cal{T}_M$. The connection
$\nabla^*$ identifies with $\nabla$ under this isomorphism.
The inverse isomorphism $(\beta^*)^{-1}$ can be described
as $X\mapsto i_X(\nabla\,Bu)$.

\smallskip

Denote $\Cal{P}^*$ the pencil of flat connections $\beta^*(\Cal{P})$.

\medskip

{\bf 4.2.1. Example.} Assume that $(M,\circ ,e, \nabla_0)$
is an $F$--manifold endowed with a compatible flat
structure and a $\nabla_0$--flat identity. Put $\Cal{F}:=\Cal{T}_M$
and construct $A$ so that $X\circ Y=i_X(A)(Y).$
Then $e$ is a primitive section which induces exactly the initial
flat structure $\nabla_0^*=\nabla_0$. 

\smallskip

In fact, let $(x^a)$ be a $\nabla_0$--flat local
coordinate system on $M$ such that $e=\partial_0$, and 
$(\partial_0, \dots ,\partial_m)$ the dual basis of flat 
vector fields. Then an easy argument shows that as a vector field,
$Be=\sum_a (x^a+c^a)\,\partial_a$ where $c^a$ are constants.

\smallskip

The following theorem borrowed from [LoMa2]
shows that a converse statement holds as well.
We will reproduce the proof here.
  
\medskip

{\bf 4.3. Theorem.} {\it Let $(M,\Cal{F},\nabla ,A,u)$
be (the data for) a pencil of flat connections on an
external bundle endowed with a primitive section.
Then $\Cal{P}^*$ is a pencil of torsionless flat connections,
and $e:=\beta^*(u)$ is an identity for one of
the associated $F$--manifold structures $\circ$.}

\medskip

{\bf Proof.} We know that 
one point $\nabla^*\in \Cal{P}^*$
is torsionless. The key remark is
that one more point is torsionless, namely $\nabla^* +A^*$
where $A^*=\beta^*(A)$. Once we have established this,
everything will follow from the results of \S 2.

\smallskip

To see this, let us work in local coordinates.
Let $\partial_a$ be a basis of $F$ considered simultaneously
as sections of $\Cal{F}$ and of $\Cal{T}_M$, such that
$\partial_0=u$. Let $(x^a)$ be dual local coordinates.
Let $A=\nabla^*B^*$ where $B^*:=\beta^*(B)$.
We represent $B^*$ as an even matrix function $(B^a_c)$
such that its action on sections of $\Cal{F}$ or $\Cal{T}_M$
is given by
$$
B\left(\sum_af^a\partial_a\right)=\sum_a\left(\sum_c B^a_cf^c\right)\partial_a.
\eqno(4.2)
$$
By construction, the map $X\mapsto i_X(\nabla^*B^*)u$
must be identical on $F$. It follows that
$B^*u$ as a vector field must be of the form
$\sum_a (x^a+c^a)\,\partial_a$ as in the example above.

\smallskip

We want to prove that
$$
i_{\partial_a}(\nabla^*B^*)(\partial_b) = (-1)^{ab}
i_{\partial_b}(\nabla^*B^*)(\partial_a)
\eqno(4.3)
$$
or equivalently
$$
\forall\,c,\ d\omega^c=0,\quad \omega^c:=\sum_adx^a\,B^c_a.
\eqno(4.4)
$$
Now, from $\nabla^*B^*\wedge \nabla^*B^* =0$
and $\nabla^*u=0$ it follows that  $\nabla^*B^*\wedge \nabla^*B^*u =0$.
Let us write this in coordinates:
$$
0=(dB^a_c)\wedge\left(\sum_a dx^a\partial_a\right)=
\sum_a\left(\sum_cdB^a_cdx^c\right)\partial_a=
\sum_a d\omega^a\partial_a.
$$
This proves (4.3) and (4.4).

\smallskip

For more details and deeper results, see [LoMa2], \S 5.

\medskip

{\bf 4.4. The tangent bundle considered as an external bundle.}
Fix now a structure $(M,\Cal{P})$ of a manifold with a unital pencil
of torsionless flat connections as in sec. 2.11. Assume moreover
that one (hence each) identity $e_A$ is flat with respect to some
$\nabla_0 \in \Cal{P}.$ Choose as origin  $\nabla_0$
and a coordinate $\lambda$ on $\Cal{P}$ so that
$(M,\Cal{P})$ determines an $F$--manifold structure $\circ$
with $\nabla_0$--flat identity $e$ and multiplication tensor $A$. 

\smallskip

In the following we will study
the family of all $F$--manifold structures that can be obtained from
this one by treating $\Cal{T}_M$ as an external bundle with the pencil $\Cal{P}$,
choosing different $\nabla \in \Cal{P}$ and different
$\nabla$--flat primitive sections, and applying Theorem 4.3.   

\smallskip

{\bf 4.4.1. Definition.} {\it An even global vector field $\varepsilon$
is called a virtual identity, if it is invertible
with respect to $\circ$ and its inverse $u:=\varepsilon^{-1}$
belongs to $Ker\,\nabla$ for some $\nabla\in\Cal{P}.$} 

\medskip

{\bf 4.4.2. Proposition.} {\it a) Inverted virtual identities 
in  $Ker\,\nabla$
are exactly primitive sections of $(\Cal{T}_M,\nabla )$
considered as an external bundle. 

\smallskip

b) The map $(\beta^*)^{-1}$ sends a local vector field $X$ to $X\circ u=X\circ\varepsilon^{-1}$.}

\smallskip

{\bf Proof.}  As we have already remarked, $(\beta^*)^{-1}$
sends a vector field $X$ to $i_X(\nabla B)(u)$.
Since $A=\nabla B$, the last expression equals $X\circ u$.
Thus $x\mapsto B(x)u$
is a local isomorphism if and only if  the 
$\circ$--multiplication by $u$ is an invertible operator,
that is, $\varepsilon$ is a virtual identity.
This completes the proof. 

\medskip

{\bf 4.5. Dubrovin's duality.} Let now $(M, \Cal{P},e,\varepsilon )$
be a flat unital pencil on $\Cal{T}_M$ endowed with a
$\nabla_0$--flat identity $e$, and a $\nabla$--flat
inverse  virtual identity $\varepsilon^{-1}$.
\medskip

{\bf 4.5.1. Theorem.} {\it a) Denote  by $*$ the new multiplication
on $\Cal{T}_M$ induced by $(\beta^*)^{-1}$ from $\circ$: 
$$
X*Y=\varepsilon^{-1}\circ X\circ Y.
\eqno(4.5)
$$
Then $(M,*,\varepsilon )$ is an $F$--manifold with identity
$\varepsilon$ (hence our term ``virtual identity'').

\smallskip

b) With the same notation, assume moreover that $\nabla_0\ne \nabla$
and that $A$ is normalized as $A=\nabla -\nabla_0.$ Then
$e$ is an Euler field of weight one for $(M,*,\varepsilon )$.}

\smallskip

{\bf Proof.} a) Clearly, $*$ is a commutative associative
$\Cal{O}_M$--bilinear multiplication. We will now exhibit
the associated  pencil  of torsionless flat connections.
We have
$$
\nabla^*_X(Y):=\varepsilon\circ \nabla_X(\varepsilon^{-1}\circ Y).
$$
More generally, the whole pencil $\Cal{P}^*
=\{\,\nabla^*_{\lambda,X} :=\nabla^*_X +\lambda X*\,\}$ is the pullback of $\Cal{P}$ 
with respect to  $Y\mapsto Y\circ u$:
$$
(\nabla_X^*+\lambda X*)Y=\varepsilon\circ \nabla_X(\varepsilon^{-1}
\circ Y) +\varepsilon \circ (\lambda \varepsilon^{-1}\circ X
\circ \varepsilon^{-1}\circ Y).
\eqno(4.6)
$$
Hence the first part follows from the Theorem 4.3.

b) We must now check the relevant versions of the Euler field relations
(3.16) and (3.17):
$$
[\varepsilon ,e] =\varepsilon ,
\eqno(4.7)
$$
$$
P^*_e(X,Y):=[e, X*Y]-[e,X]*Y-X*[e,Y]=X*Y.
\eqno(4.8)
$$
In fact, we have $[\varepsilon ,e]=-\nabla_{0e}\varepsilon$
because $\nabla_0e=0$. Moreover, $\nabla_0=\nabla -A$ and $\nabla\varepsilon =0$,
hence $\nabla_{0e}\varepsilon =-e\circ\varepsilon =-\varepsilon$
which shows (4.7).

\smallskip

Furthermore,
$$
P^*_e(X,Y)=[e, \varepsilon^{-1}\circ X\circ Y]-\varepsilon^{-1}\circ [e,X]
\circ Y-\varepsilon^{-1}\circ X \circ [e,Y].
$$
Using the fact that $\roman{ad}\,e$ is a $\circ$--derivation (cf. (2.10)),
we can rewrite the first summand in the right hand side. After cancelations
we get
$$
P^*_e(X,Y)= [e,\varepsilon^{-1}]\circ X\circ Y.
$$
Finally, 
$$
[e,\varepsilon^{-1}]=-\varepsilon^{-2}\circ [e,\varepsilon ]=\varepsilon^{-1}
$$
in view of (4.7), which completes the check of (4.8).

\smallskip

The last property to be checked is $[e, \roman{Ker}\,\nabla^*_0]\subset \roman{Ker}\,\nabla^*_0$ (notice that $\varepsilon$
is $\nabla^*_0$--flat, because $e$ is $\nabla_0$--flat). But $\roman{Ker}\,\nabla^*_0=\varepsilon\circ\roman{Ker}\,\nabla_0$,
and for a $\nabla_0$--flat $X$ we have $[e,\varepsilon\circ X]=
[e,\varepsilon ]\circ X=-\varepsilon\circ X .$

\medskip

{\bf 4.5.2. Comments.} a) The relationship
between $(M,\circ , e)$ and $(M,*,\varepsilon )$
is almost symmetric, but not quite. To explain this, we
start with a part that admits a  straightforward check:
$$
X\circ Y=e^{*-1}* X*Y
\eqno(4.5)
$$
which is the same as (4.1) with the roles of $(\circ ,e)$
and $(*,\varepsilon )$ inverted. Here $e^{*-1}$ 
denotes the solution $v$ to 
the equation $v*e=\varepsilon$, that is  
$\varepsilon^{-1}\circ v\circ e=
\varepsilon$, therefore $v=e^{*-1}=\varepsilon^{\circ 2}.$
After this remark one sees that (4.5) follows from (4.1).

\smallskip

Slightly more generally, one easily checks that 
the map inverse to
$X\mapsto \varepsilon^{-1}\circ X$ reads
$Y\mapsto e^{*-1}*Y$ as expected.

\smallskip 

If the  symmetry were perfect, we would now expect $\varepsilon$ 
to be an Euler field of
weight one for $(M,\circ ,e)$, however, this contradicts 
(4.3)!

\smallskip

The reason of this is that the vector field
$e$ is not a virtual identity for $(M,*,\varepsilon )$:
$e^{*-1}=\varepsilon^2$ is not flat with respect to any
$\nabla^*\in \Cal{P}^*$, so that if we start with $(M,*,\varepsilon)$
and construct $\circ$ via (4.5), we cannot apply
Theorem 4.5.1 anymore.

\smallskip

This remark suggests two ways of extending our construction.

\smallskip

First, one can iterate it by applying it
to $(M,*,\varepsilon, \eta )$ where $\eta$
is an arbitrary  virtual identity thus getting new structures
of an $F$--manifold.

\smallskip

Second, and more interesting, one can try to extend the definition of
a virtual identity: say, call {\it an eventual identity}
for $(M,\circ ,e)$ any invertible vector field
$\varepsilon$ such that (4.1) defines a structure of
an $F$--manifold. We have seen that virtual
identities are eventual ones, but not netcessarily vice versa.
Can one give an independent characterization of
eventual identities?

\medskip

b) If we take a pair consisting
of a Frobenius manifold and the dual almost Frobenius
manifold in the sense of [Du2], and retain only
their $F$--manifold structures with the relevant
flat structures, we will get a pair like ours
$(M,\circ , e)$ and $(M,*,\varepsilon )$.

\smallskip

A reader willing to compare our work with [Du2],
should look at this point at Dubrovin's formulas (3.1) and (3.26)
and replace Dubrovin's $*,\cdot , E, e, u, v$ respectively by our
$\circ,*,e, \varepsilon, X, Y$.

 \smallskip

Our extension shows that the use of a metric
is superfluous, so that the term
``duality'' is not quite justified in this context,
because the construction is more like ``twisting''.
Moreover, in the realm of $F$--manifolds
with compatible flat structure the construction becomes
somewhat more transparent and natural.

\bigskip

\centerline{\bf \S 5. Formal flat $F$--manifolds and toric compactifications}

\medskip

{\bf 5.1. Notation.} Let $T$ be a finite--dimensional linear superspace
over $\bold{C}$ (any characteristic zero ground field 
or a local Artin superalgebra $k$ over it will do as well).
Consider a basis $(\Delta_a\,|\,a\in I)$ of $T$ and a
system of linear coordinates $(x^a)$ on
$F$ and denote by $M$ the formal completion of $T$ at 0, that is,
the formal spectrum of $\Cal{R}:=\bold{C}[[x^a]].$ 
The space of vector fields $T_M$ on $M$ can be
canonically identified with $\Cal{R}\otimes T$:
$\Delta_a\mapsto \partial_a.$
Let $\nabla$ be the torsionless free connection 
on $T_M$ with kernel $T$. 

\smallskip

According to the results of \S 2, classification of
all formal $F$--manifold structures on $M$
compatible with $\nabla$ is equivalent
to the classification  of the solutions
to the matrix differential equation
$$
\nabla B\wedge \nabla B=0,\quad B\in m\otimes_{\bold{C}} \roman{End}\,T,
\eqno(5.1)
$$
where $m$ is the maximal ideal of $\Cal{R}$. 

\smallskip

Namely, given such a solution, we put
$$
X\circ Y:=i_X(\nabla B)(Y)=(XB)Y
\eqno(5.2)
$$
where $X$ acts as derivation on the  entries of $B$. 
 
\smallskip

The pencil of torsionless flat connections associated with
$(M,\circ ,\nabla )$ is $\nabla +\lambda \nabla B$.

\smallskip

This remark allows us to apply to this situation a
result of [LoMa1] which states essentially
that the classification of solutions to (5.1) 
reduces to a problem in the representation theory of
an associative algebra. We will summarize
this result below, refering to [LoMa1]
for proofs, and for more general statements
about connections on an external bundle.

\smallskip

Notice that there exists another description of a
formal $F$--manifold with compatible flat structure on $M$:
it is the same as the structure of an operadic
algebra over the {\it oriented} homology operad
$(H_{*,or}(\overline{M}_{0,n+1})$. Orientation means that
the zeroth structure section of each 
universal curve $C_{0,n+1}\to \overline{M}_{0,n+1}$
is marked, the underlying combinatorial
formalism involves rooted trees with root corresponding to
the $0$--th section, and
grafting is allowed only root--to--nonroot.
The $n$--th component of the operad $H_{*,or}(\overline{M}_{0,n+1})$
is acted upon by the $n$--th symmetric group
$\bold{S}_n$, rather than $\bold{S}_{n+1}$.

\smallskip

For details of this operadic description and its
generalizations,  see [LoMa2].

\medskip

{\bf 5.2. $\bold{S}$--algebras.} Consider 
a graded associative $\bold{C}$--algebra $V=\oplus_{n=1}^{\infty}V_n$
(without identity) in
the category of vector superspaces.
I will call it  {\it an $\bold{S}$--algebra}, if for each $n$,
an action of the symmetric group $\bold{S}_n$ on $V_n$
is given such that the multiplication map $V_m\otimes V_n\to V_{m+n}$
is compatible with the action of $\bold{S}_m\times\bold{S}_n$
embedded into $\bold{S}_{m+n}.$ Example: 
the tensor algebra
(without the rank zero part) of a superspace $T$ is an  
$\bold{S}$--algebra.

\smallskip

In any  $\bold{S}$--algebra $V$ the sum of subspaces
$J_n$ spanned by $(1-s)v,\,s\in\bold{S}_n,\,v\in V_n,$
is a double--sided ideal in $V.$ Hence the sum of
the coinvariant spaces $V_{\bold{S}}:=\oplus_n V_{\bold{S}_n} =V_n/J_n$ is a graded ring. I will denote it $V_{\bold{S}}$.

\smallskip

If $V$, $W$ are two $\bold{S}$--algebras, then the diagonal
part of their tensor product $\oplus_{n=1}^{\infty}V_n\otimes W_n$
is an  $\bold{S}$--algebra as well.

\medskip

{\bf 5.3. Algebra $H_*$.} We will now define an $\bold{S}$--algebra
whose $n$--th component is the homology of the toric
variety associated with the $n$--th permutohedron.

\smallskip

{\bf 5.3.1. Permutohedral fans.} Let $B$ be a finite
set. A partition $\sigma$ of $B$ is defined as {\it a totally
ordered set of pairwise disjoint non--empty subsets of $B$ whose union is
$B$.}
If a partition consists of $N$ subsets, it is called an
$N$--partition. If its components are denoted $\sigma_1,\dots ,\sigma_N$,
this means that they are listed in their structure order.

\smallskip

Let $\tau$ be an $N+1$--partition of $B$. If $N\ge 1,$ it determines
a well ordered family of $N$ 2--partitions $\sigma^{(a)}$:  
$$
\sigma^{(a)}_1:=\tau_1\cup\dots\cup\tau_{a},\
\sigma^{(a)}_2:=\tau_{a+1}\cup\dots\cup\tau_{N},\ a=1,\dots ,N\, .
$$
Call a sequence of $N$ 2--partitions $(\sigma^{(i)})$ 
{\it good} if it can be obtained by such a construction.

\smallskip

Put $N_B := \bold{Z}^B/\bold{Z}$, the latter subgroup being embedded
diagonally. Similarly, $N_B\otimes{\bold{R}}=\bold{R}^B/\bold{R}$.
The vectors in this space (resp. lattice) will be written as functions
$B\to \bold{R}$ (resp. $B\to \bold{Z}$) considered modulo
constant functions. For a subset $\beta\subset B$, let $\chi_{\beta}$
be the function equal 1 on $\beta$ and 0 elsewhere.

\smallskip 

The fan $\Phi_B$ in $N_B\otimes{\bold{R}}$, consists
of certain $l$--dimensional cones $C(\tau )$ labeled by 
$(l+1)$--partitions $\tau$ of $B$.

\smallskip

If $\tau$ is the trivial 1--partition, $C(\tau )=\{0\}$.

\smallskip

If $\sigma$ is a 2--partition, $C(\sigma )$ is generated
by $\chi_{\sigma_1}$, or, equivalently, $-\chi_{\sigma_2}$,
modulo constants.

\smallskip

Generally, let $\tau$ be an $(l+1)$--partition, and 
$\sigma^{(i)},\,i=1,\dots,l$,
the respective good family of 2--partitions. Then
$C(\tau )$ is defined as a cone generated by all $C(\sigma^{(i)})$.

\medskip

{\bf 5.3.2. Permutohedral toric varieties.} Denote by $L_B$ 
the toric variety which is the compactification
of the torus $(\bold{G}_m)^B/\bold{G}_m$
associated with the fan $\Phi_B$. One can prove that it is smooth and proper;
its dimension is $\roman{card}\,B -1.$ By construction,
the permutation group of $B$ acts upon it.

\medskip

{\bf 5.3.3. The $\bold{S}$--algebra $H_*$.}
By definition, the  $n$--th component of $H_*$ is the homology
space $H_{*n}:=H_*(L_n)$ where we write now $L_n$ for
$L_{\{1,\dots ,n\}}$. The whole homology space
is spanned by certain cycles $\mu (\tau )$
indexed, as well as cones of $\Phi_n$, by partitions $\tau$ of
$\{1,\dots ,n\}$: this is a part of the general
theory of toric compactifications. 
 
\smallskip

Define now a multiplication $H_*(L_m)\times H_*(L_n)\to H_*(L_{m+n}):$
if $\tau^{(1)}$ (resp. $\tau^{(2)}$)
is a partition of $\{1,\dots ,m\}$ (resp. of $\{1,\dots ,n\}$),
then
$$
\mu (\tau^{(1)})\,\mu (\tau^{(2)}) := \mu (\tau^{(1)}\cup \tau^{(2)})
$$
where the concatenated partition of 
$\{1,\dots ,m,\, m+1,\dots ,m+n\}$ is defined in an obvious way,
shifting all the components of $\tau^{(2)}$ by $m$.
Of course, this map has a geometric origin coming from
certain boundary morphisms $L_m\times L_n\to L_{m+n}.$

\medskip

{\bf 5.4. Algebra $H_*T$,\, its representations
and matrix correlators.} Let now $T$ be a linear superspace considered 
in 5.1. Define  $H_*T$ as the algebra of coinvariants of the 
diagonal tensor product
$$
H_{*}T:= \left(\oplus_{n=1}^{\infty} H_{*n}\otimes T^{\otimes n}
\right)_{\bold{S}} .
$$
Consider a linear representation $\rho:\,H_*T\to \roman{End}\,T$
and define {\it the matrix correlators}
of $\rho$ as the following family of endomorphisms of $T$:
$$
\tau^{(n)}\langle \Delta_{a_1}\dots\Delta_{a_n}\rangle_{\rho} :=
\rho(\mu (\tau^{(n)})\otimes \Delta_{a_1}\otimes\dots\otimes\Delta_{a_n}) .
$$
Here $\tau^{(n)}$ runs over all partitions of $\{1,\dots ,n\}$
whereas $(a_1,\dots ,a_n)$ runs over all maps $\{1,\dots ,n\}\mapsto I:\,
i\to a_i.$

\smallskip

{\it Top matrix
correlators of $\rho$} constitute the subfamily of correlators corresponding
to the identical partitions $\varepsilon^{(n)}$ of $\{1,\dots ,n\}$:
$$
\langle \Delta_{a_1}\dots\Delta_{a_n}\rangle_{\rho} :=
\varepsilon^{(n)}\langle \Delta_{a_1}\dots\Delta_{a_n}\rangle_{\rho}\, .
$$
The relevant cycle $\mu (\varepsilon^{(n)})$ is the fundamental
cycle of $L_n$.
 
\medskip

{\bf 5.5. Main theorem.} We can now state the basic correspondence
between the solutions to (5.1) and representations ${\rho}$.

\smallskip

Given ${\rho}$, construct the series
$$
B_{\rho} =
\sum_{n=1}^{\infty}\sum_{(a_1,\dots ,a_n)}
\frac{x^{a_n}\dots x^{a_1}}{n!}\,
\langle \Delta_{a_1} \dots \Delta_{a_n}\rangle_{\rho} \in \bold{C}[[x]]\otimes \roman{End}\,F.
\eqno(5.3)
$$

\medskip

{\bf 5.5.1. Theorem.} {\it a) We have
$$
\nabla B_{\rho}\wedge \nabla B_{\rho}=0.
\eqno(5.4)
$$

\smallskip

b) Conversely, let $\Delta (a_1,\dots ,a_n)\in\roman{End}\,F$ be a
family of linear operators defined for all $n\ge 1$ and all
maps $\{1,\dots ,n\}\to I:\,i\mapsto a_i$. Assume that the parity
of $\Delta (a_1,\dots ,a_n)$ coincides with the sum of the
parities of $\Delta_{a_i}$ and that 
$\Delta (a_1,\dots ,a_n)$ is (super)symmetric with respect
to permutations of $a_i$'s.
Finally, assume that the formal series 
$$
B =
\sum_{n=1}^{\infty}\sum_{(a_1,\dots ,a_n)}
\frac{x^{a_n}\dots x^{a_1}}{n!}\,
\Delta ({a_1}, \dots ,{a_n}) \in \bold{C}[[x]]\otimes \roman{End}\,F
\eqno(5.5)
$$
satisfies the equations (5.1). Then there exists a well defined
representation ${\rho}:\,H_*T\to \roman{End}\,F$ such that
$\Delta ({a_1}, \dots ,{a_n})$ are the top correlators
$\langle \Delta_{a_1} \dots \Delta_{a_n}\rangle_{\rho}$ of this
representation.}

\bigskip

\centerline{\bf References}

\medskip

[Du1] B.~Dubrovin. {\it Geometry of 2D topological field theory.}
In: Springer Lecture Notes in Math. 1620 (1996), 120--348.

\smallskip

[Du2] B.~Dubrovin. {\it On almost duality for Frobenius manifolds.}
Preprint \newline math.DG/0307374.

\smallskip

[Gi] A.~Givental. {\it On the WDVV--equation in quantum $K$--theory.}
Fulton's Festschrift, Mich. Math. Journ., 48 (2000), 295--304.

\smallskip

[He] C.~Hertling. {\it Frobenius manifolds and moduli spaces for singularities.}
Cambridge University Press, 2002.

\smallskip

[HeMa] C.~Hertling, Yu.~Manin {\it Weak Frobenius manifolds.}
Int. Math. Res. Notices, 6 (1999), 277--286. Preprint
math.QA/9810132.

\smallskip

[Lee] Y.~P.~Lee. {\it Quantum $K$--theory I: foundations.}
Preprint math.AG/0105014.

\smallskip

[LoMa1] A.~Losev, Yu.~Manin. {\it New moduli spaces of pointed curves and pencils of flat connections.} Fulton's Festschrift,
Michigan Journ. of Math., 48 (2000), 443--472. Preprint math.AG/0001003

\smallskip

[LoMa2] A.~Losev, Yu.~Manin. {\it Extended modular operad.} Preprint 
\newline math.AG/0301003.

\smallskip

[Ma] Yu.~Manin. {\it Frobenius manifolds, quantum cohomology, and moduli spaces.} AMS Colloquium Publ. 47, Providence RI, 1999, 303 pp.

\smallskip

[Me1] S.~Merkulov. {\it Operads, deformation theory and $F$--manifolds.}
Preprint math.AG/0210478.

\smallskip

[Me2] S.~Merkulov. {\it PROP profile of Poisson geometry.}
Preprint math.AG/0401034.

\enddocument